\documentclass[a4wide, 11pt]{article}

\usepackage[french]{babel}
\usepackage[utf8]{inputenc}
\usepackage{enumitem}
\usepackage{amsmath}
\usepackage{amssymb}
\usepackage{titlesec}
\usepackage{pdfpages}
\usepackage[colorlinks]{hyperref}

\newcommand{\Nm}{\ensuremath{\mathbb{N}}}

\newcommand{\Cm}{\ensuremath{\mathbb{C}}}

\newcommand{\vs}{\vspace{.2cm}}

\newtheorem{thm}{Theorem}

\def\qed {\mbox{}\hfill {\small \fbox{}} \\}
\def\lto{\longrightarrow}
\def\lmto{\longmapsto}

\def\leq{\leqslant}
\def\geq{\geqslant}

\titleformat{\subsubsection}[runin]
{\normalfont\large\bfseries}{\thesubsubsection}{1em}{}
\counterwithout{subsubsection}{subsection}

\renewcommand{\thesubsubsection}{(\arabic{subsubsection})}

\title{The Siegel-Bruno linearization Theorem}
\author{Patrick Bernard}
\date{}

\begin{document}

\maketitle
The purpose of this paper is to provide a short and self-contained account on  Siegel's Theorem, as improved by Bruno, which states that a holomorphic map $f$ of $\Cm$ which fixes $0$ can be locally linearized, under certain conditions on the multiplier $\lambda:=f'(0)$.

\begin{thm}[Koenigs, Siegel, Bruno]
	Let $f:U\subset \Cm\lto \Cm $ be a holomorphic map defined on an open set $U$ containing the origin, such that $f(0)=0$. Suppose furthermore that $|\lambda|\not \in \{0,1\}$ (this is the hyperbolic case),
	or that $\lambda\in B\subset S^1$ (this is the elliptic case), where $B$ will be described below. Then there exists a unique holomorphic diffeomorphism  $h$, defined in the neighborhood of $0$, such that $h'(0)=1$ and such that $h^{-1}\circ f\circ h=\lambda I$ in the neighborhood of $0$. 
\end{thm}

The hyperbolic case, attributed to Koenigs \cite{K}, dates from the 19th century. The elliptic case is more difficult, Siegel gave in 1942 a  proof for a set $D\subset S^1$ of total Lebesgue measures, the Diophantines. Bruno then obtained the result for a larger set $B\supset D$ of multipliers. Yoccoz finally showed in \cite{Yo} that the set $B$ is optimal, the map $f(z)=\lambda z+z^2$ is not linearizable if $\lambda \in S^1-B$.

To define the set $B$, let us introduce the small divisors
$$
\omega_n=|\lambda^n-1|, \quad \Omega_n=\min _{1\leq \ell\leq n}\omega_{\ell}.
$$ 
An important difference between the hyperbolic case and the elliptic case is that these small divisors are bounded from below by a strictly positive real in the hyperbolic case. On the contrary, in the elliptic case, the sequence $\Omega_n$ converges to $0$. The conditions to obtain the linearization express that this convergence is not too fast (in particular, this sequence should be positive, and thus $\lambda$ should not be a root of the unity). More precisely the set  $B\subset S^1$ of Bruno multipliers is defined by :
$$\lambda \in B \quad \Longleftrightarrow \quad \sum _{k\geq 1}2^{-k}\ln \Omega_{2^k}^{-1}<\infty.
$$
The multiplier $\lambda$ is said to be Diophantine if there exist $a>0$ and $b>0$ such that 
$$
\Omega_n\geq an^{-b},
$$
it is immediate then that $\lambda\in B$. 

The classical methods to prove the theorem consist either in studying the formal series defining the conjugacy (which is well defined as soon as $\lambda$ is not a root of unity) or in studying an iterative approximation of the conjugacy. The first method is  used in Siegel's original article, \cite{S}, it is also exposed (for Bruno multipliers) by Bruno in \cite{B}, \S 5. The second method is used by Rüssmann  in \cite{R}, and is also exposed (for Diophantine multipliers) in the books \cite{CG}, II.6, or \cite{A} \S 28. Another proof by renormalization is given by Yoccoz in \cite{Yo}.
We present here an intermediate method, taken from Bruno's article \cite{B}, \S 4 (where it is used in the context of differential equations in higher dimension) which consists in studying the series  by an iterative method. The present paper does not introduce novelties compared to Bruno's iterative procedure, but rather attempts to give a pedagogical self-contained account of this method by presenting it in a simple situation. Several improvements have been obtained since  the paper of Bruno, it is not our purpose here to review this very rich literature.

\subsection{Notations}
Let $\Cm[[z]]$ be the set of complex power series $f=\sum(f)_kz^k$. An element of $\Cm[[z]]$ is thus a complex sequence $(f)_k, k\in \Nm$. 
It is a $\Cm$-algebra for  the Cauchy product $f\cdot g$ given by
$$
(f\cdot g)_k:=\sum_{i=0}^k (f)_i (g)_{k-i}.
$$

For $z\in \Cm$ and $f\in\Cm[[z]]$, we denote  by$f(z)$ the value of the series 
$$
\sum _{k=0} ^{\infty} (f)_k z^k
$$
when it converges. We denote by $\rho(f)\in [0, \infty]$  the radius of convergence of the series $f$.

Let $O_k\subset \Cm[[z]]$ be the space of power series
whose $k$ first coefficients are zero, i.e. series of the form $(f)_kz^k+(f)_{k+1} z^{k+1}+\cdots$.
For $f\in \Cm[[z]]$ and $d\in \Nm$, we denote
$$
[f]_d:=(f)_0+(f)_1z+\cdots +(f)_d z^d
$$
the polynomial obtained by truncating $f$ to order $d$. It is an element of $\Cm[[z]]$, and $f-[f]_d\in O_{d+1}$.

We will consider sequences $f_n$ of elements of $\Cm[[z]]$, and we will say that a sequence converges strongly to $f$
if each of the sequences $n\lmto (f_n)_k$ stabilizes at the value $(f)_k$ for large $n$.

When $f :U\subset \Cm\lto \Cm$ is an analytic function, we still denote by $f$ the corresponding power series.
For example, we denote $1/(1-z)$ the series $1+z+z^2 +\cdots$. We denote by $I$, or $z$, the identity map of the complex plane as well as the associated power series.

Given a power series $f$, we denote by  $\hat f$  the power series whose coefficients are the moduli of the coefficients of $f$. 
For $f, g\in \Cm[[z]]$, it is easy to verify that
$$
\widehat {f\cdot g} (r)\leq \hat f(r)\hat g(r).
$$

\subsection{Composition}

If $f$ is any power series and $g$ is a power series with no constant term (i.e. $g\in O_1$), we define as usual the composition $f\circ g$ by 
$$
(f\circ g)_m=\sum_{k=0}^{m} (f)_k (g^k)_m,
$$
where $(g^k)_m$ is the coefficient of degree $m$ of the product $g^k=g\cdot g \cdots g$ ($k$ factors) for $k\geq 1$, and $g^0=1$.

\subsubsection{}
The following formal calculation
$$f\circ g(z)=\sum_m \sum _{k=0}^m (f)_k(g^k)_m z^m= \sum _k(f)_k \sum_{m\geq k}(g^k)_mz^m=\sum _k(f)_k g(z)^k=f(g(z))
$$
justifies the notation. It is correct when the two-index family $(f)_k(g^k)_mz^m$ is summable, this is the case when $\hat f(\hat g(|z|))<\infty$.
In summary:

\vs
\noindent
{\bf Property }
\begin{itshape}
	If $z\in \Cm$ is such that $\hat f(\hat g(|z|))<\infty$, then $g$ converges at $z$, $f$ converges at $g(z)$, 
	$f\circ g$ converges at $z$, 
	and
	$f(g(z))=f\circ g(z)$.
\end{itshape}
\vs 

This implies in particular that $f\circ g$ has  positive radius if $f$ and $g$ have  positive radius.
\subsubsection{}
The next properties are easy to verify:
\begin{itemize}
	\item $f\in O_n, g\in O_k\Rightarrow f\cdot g\in O_{k+n},$
	\item $[f\cdot g]_d=[[f_d]\cdot [g]_d]_d,$
	\item $f\in O_n, g\in O_k, k\geq 1\Rightarrow f\circ g\in O_{kn},$
	\item $f\in O_n, g\in O_1, h\in O_k \Rightarrow f\circ (g+h)-f\circ g \in O_{n+k-1}$.
	\item $[f\circ g]_d=[[f_d]\circ [g_d]]_d.$
\end{itemize}
From this last property, we deduce the associativity of the composition:
$$
(f\circ g)\circ h=f\circ (g\circ h).
$$
Indeed, this property is satisfied by the truncated series $[f]_d, [g]_d, [h]_d$, which are polynomials
(for which the composition is the usual composition of polynomials), and so the truncations at all orders of the two members of the equality are equal.

\subsection{Inverse series}

\subsubsection{}\label{rec}
{\bf Proposition}
\begin{itshape}
	Let $f=\lambda I+F$, $\lambda \neq 0$, $F\in O_2$. Then there exists a unique series $G\in O_2$ such that
	$(\lambda I+F)\circ(\lambda^{-1} I+G)=I$ and
	$(\lambda^{-1} I+G)\circ(\lambda I+F)=I$.
	
	Moreover, if there exist $r>0$ and $\alpha\in ]0,1[$ such that $\hat F(r)\leq |\lambda| \alpha r$, then 
	$$\hat G(|\lambda|(1-\alpha)r)\leq \alpha r.
	$$
	\end{itshape}

We denote by  $f^{-1}$  the inverse series $\lambda^{-1}I+G$.
If $F$ has a positive radius of convergence, then $\hat F'(0)=0$ so the condition $\hat F(r)\leq |\lambda|\alpha r$
is satisfied for $r$ small, and so $f^{-1}$ has a  positive radius of convergence.

\subsubsection{Proof. }
We denote $\mu:=\lambda^{-1}$.
The equation $(\lambda I+F)\circ (\mu I+G)=I$ can be rewritten
$$
G=-\mu F\circ (\mu I+G).
$$
As usual in a fixed point problem, we consider the sequence 
$G_k $ of elements of $O_2$ defined by recurrence by 
$G_{k+1}= -\mu F\circ(\mu I+G_k) $, the first term $G_1$ being any element of $O_2$.

We will show by recurrence that this sequence converges strongly to a limit $G$ which does not depend on $G_1$, and which satisfies the inversion equation.
More precisely, we show the recurrence hypothesis: $G_k\in G_{k-1}+O_k$ and $[G_k]_k$ does not depend on $G_1$.

Denoting $O_k$ any element of $O_k$, if $G_k=G_{k-1}+O_{k},$ then, as $F\in O_2$,
$$G_{k+1}=- \mu F\circ (\mu I+G_{k-1}+O_{k})=- \mu F\circ (\mu I+G_{k-1})+O_{k+1}=G_k+O_{k+1}.
$$
In a similar way,
$$G_{k+1}=- \mu F\circ (\mu I+[G_{k}]_k+O_{k+1})=- \mu F\circ (\mu I+[G_k]_k)+O_{k+2},
$$
so $[G_{k+1}]_{k+1}$ depends only on $[G_k]_k$, and therefore does not depend on $G_1$ by the recurrence hypothesis.

This implies that the sequence $G_k$ converges strongly to a limit $G$, characterized by $[G]_k=[G_{k}]_k$.
Then for all $k$ we have
\begin{align*} [G]_k&=[G_{k+1}]_k=[-\mu F\circ (\mu I+G_{k})]_k=
	[-\mu F\circ (\mu I+G+O_{k+1})]_k=[-\mu F\circ (\mu I+G)]_k,
\end{align*}
so the equality $G=-\mu F\circ (\mu I+G)$ is satisfied by the limit.

We can present the above in a slightly different way. 
Since only the truncation $[G_k]_k$ matters, we could consider the sequence of polynomials 
$G_k:=[G_k]_k$, defined by recurrence by $\tilde G_1=0$ and 
$$
\tilde G_k=[-\mu F\circ (\mu I+\tilde G_{k-1})]_k.
$$
We can verify as above by recurrence that $\tilde G_k\in \tilde G_{k-1}+O_k$, 
that is to say that the passage from $\tilde G_{k-1}$ to $\tilde G_k$ just consists in adding a term of order $k$, which is given by the recurrence relation. This is the classical proof  of the existence and uniqueness of the formal series we are looking for.

To prove that the right inverse  is equal to the left inverse, we can consider the right inverse  $(\lambda I+H)$ of $(\mu I+G.)$
Then, 
$\lambda I+F=(\lambda I+F)\circ (\mu I+G)\circ (\lambda I+H)=\lambda I+H,
$
so $H=F$.

We finally show by recurrence that $\hat G_k(|\lambda|(1-\alpha)r)\leq \alpha r$.
Assuming this recurrence hypothesis, we have 
$$
\hat G_{k+1}(|\lambda|(1-\alpha)r)\leq |\mu| \hat F ( (1-\alpha)r+\alpha r)\leq|\mu| \hat F(r) \leq \alpha r.
$$
This implies in particular that $[\hat G]_k(|\lambda|(1-\alpha)r)\leq \alpha r$ for all $k$ and thus that $\hat G(|\lambda|(1-\alpha)r)\leq \alpha r$.
\qed

\subsection{Formal Linearization}\label{form}
	
	\subsubsection{Proposition}\textit{
		Let $f=\lambda z+F$ be a formal series, with $F\in O_2$. 
		If $\lambda \neq 0$ is not a root of unity, there exists a unique formal series $h$ of the form $h=I+H, H\in O_2$, such that
		$h^{-1}\circ f\circ h=\lambda I$.}

	\subsubsection{Proof. } The conjugacy equation is written again 
	$$
	H\circ (\lambda I) - \lambda H = F\circ (I+H).
	$$ 
	We notice that the linear operator $L_{\lambda }:H\mapsto H\circ (\lambda I) -\lambda H$
	is diagonal,
	$$
	(L_{\lambda}H)_m:=(H\circ (\lambda I) -\lambda H)_m=(\lambda^m-\lambda)(H)_m.
	$$
	If $\lambda$ is not a root of the unit, the coefficients of $L_{\lambda }$ are non-zero and
	we can treat the equation 
	$$
	H=L_{\lambda}^{-1}\big(F\circ (I+H)\big)
	$$
	exactly as the inversion equation. We define the sequence 
	$H_k$ by $H_1\in O_2$ and $H_{k+1}=L_{\lambda}^{-1}(F\circ (I+H_k))$.
	We verify as above by recurrence that $H_{k}-H_{k-1}\in O_k$, and thus that $H_k$ converges strongly to a series $H$, which
	satisfies the conjugacy equation.
	Indeed, if $H_k=H_{k-1}+O_k$ then, as $F\in O_2$,
	\begin{align*}
		H_{k+1}&=L_{\lambda}^{-1} (F\circ ( I+H_{k-1}+O_k))\\
		&= L_{\lambda}^{-1} ( F\circ ( I+H_{k-1})+O_{k+1})
		=H_k+O_{k+1}.
	\end{align*}
	Moreover, the limit $H$ does not depend on $H_1$, and it is therefore the unique solution of the equation.
	
\subsection{Linearization, hyperbolic case}
\subsubsection{} 

We prove here the hyperbolic case of the theorem: If $\rho(f)>0$ and if $|\lambda|\not \in \{0,1\}$, then $\rho(h)>0$.

We set $\omega=\inf_{m\geq 2}|\lambda^m-\lambda|$. The specificity of the hyperbolic case is that $\omega >0$, the 
operator $L_{\lambda}^{-1}$ is therefore bounded, and
we can study the convergence of the conjugacy $H$ exactly as  the inverse $G$. We obtain, more precisely:

\vs
\noindent
{\bf Proposition}
\begin{itshape}
If $\hat F(\omega^2 r)\leq \alpha \omega  r$ for some $ \alpha\in ]0,\omega^2 [$, then 
$\hat H((\omega^2 -\alpha)r)\leq \alpha r$.
\end{itshape}
\vs

As above, it is sufficient to show by recurrence that $\hat H_k((\omega^2-\alpha)r)\leq \alpha r$, which follows from the calculation
\begin{align*}
\hat H_{k+1}((\omega^2-\alpha )r)
&\leq \omega^{-1}\hat F\big((\omega^2-\alpha )r+\hat H_k ((\omega^2-\alpha )r)\big)
\leq \omega^{-1}\hat F((\omega^2-\alpha )r+\alpha r)\\
&\leq    \omega^{-1}\hat F(\omega^2 r)\leq \omega^{-1} \alpha \omega r=\alpha r.
\end{align*}
Of course, if $\rho(F)>0$, then for any $\alpha\in ]0, \omega^2[$ there exists $r>0$ such that 
$\hat F(\omega ^2r)\leq \alpha \omega r$, because $\hat F'(0)=0$.
We deduce that $\rho(H)\geq (\omega^2-\alpha)r>0$.

\subsection{Linearization, elliptic case}
We now study the linearization problem in the case $|\lambda|=1$.

\subsubsection{}
We first describe  another iterative construction of the conjugacy, which will allow a better convergence study.
We start  as earlier by posing $P=L_{\lambda}^{-1}F$.
We then check that 
$$
(I+P)^{-1}\circ (\lambda I+F)\circ (I+P)\in \lambda I+O_3.
$$
But there is more: if we already have $F\in O_{m+1}$, $m\geq 1$, then
$$
(I+P)^{-1}\circ (\lambda I+F)\circ (I+P)\in \lambda I+O_{2m+1}.
$$
To verify this, we denote $I+R:= (I+P)^{-1}$. We have 
$I=(I+R)\circ (I+P)=I+P+R\circ(I+P)=I+P+R+O_{2m+1},$ 
so 
$R+P\in O_{2m+1}$. Then, 
\begin{align*}
	(I+P)^{-1}\circ (I+F)\circ (\lambda I+P)&= (I-P)\circ (\lambda I+F)\circ (I+P)+O_{2m+1}\\
	&= \lambda I -P\circ (\lambda I)+\lambda P +F+O_{2m+1}\\
	&=\lambda I+F-L_{\lambda} P+O_{2m+1}=\lambda I+O_{2m+1}.
\end{align*}
The same calculation shows that we can replace $P$ by any series equal to $L_{\lambda}^{-1}F$ modulo $O_{2m+1}$, in particular by 
$[L_{\lambda}^{-1}F]_{2m}$.

In view of these remarks, an iterative procedure appears natural:
We pose $F_0=F$, $P_0=[L_{\lambda}^{-1}F]_{2}$, so that 
$$
F_1:= (I+P_0)^{-1}\circ (\lambda I+F_0)\circ (I+P_0)\in \lambda I+ O_{3}.
$$
Then, we apply the same procedure to the map $\lambda I+F_1$, exploiting that $F_1\in O_3$, i.e. we take $P_1:=[L_{\lambda}^{-1}F_1]_{4}$,
so that 
$$
F_2:= (I+P_1)^{-1}\circ (\lambda I+F_1)\circ (I+P_1)-\lambda I\in O_{5},
$$
and so on.
We thus define the sequences
$$
P_{k}=[L_{\lambda}^{-1}F_{k}]_{2^{k+1}}, \quad
F_{k+1}=(I+P_k)^{-1}\circ (\lambda I+ F_{k})\circ (I+P_k)-\lambda I,
$$
and verify iteratively with the help of the previous remarks that 
$$
F_k \in O_{1+2^k}, \quad P_k\in O_{1+2^k}.
$$
Setting 
$$
h_k:=(I+P_0)\circ (I+P_{1})\circ \cdots \circ (I+P_{k-1}),
$$
we obtain
$$
\lambda I+F_{k+1}=(I+P_k)^{-1}\circ (\lambda I+F_k)\circ (I+P_k)=h_k^{-1}\circ (\lambda I+F)\circ h_k.
$$
As $h_{k+1}=h_k\circ (I+P_{k})\in h_k+O_{1+2^k}$, the sequence $h_k$ converges to a limit $h$, which satisfies $\lambda I=h^{-1}\circ (\lambda I+F)\circ h$
and which is therefore the formal conjugacy. The convergence is much faster than the previous construction, since $h_k$ 
is equal to $h$ at order $2^k$ (against $k$ for the first construction). This is called quadratic convergence.

\subsubsection{}
We will study the convergence of $h$ by an inductive procedure.
We assume that $F$ is convergent (i.e. $\rho(F)>0$), and we fix, once and for all, a real $r_0>0$ such that 
$\hat F(r_0)\leq r_0$. Such a real exists because $\hat F'(0)=0$. We have $\hat F(r)\leq r$ for all $r\in [0, r_0]$.
We will prove that 
$$
\widehat{F_k}(r_k)\leq r_k, \quad \widehat{P_k}(r_{k+1})\leq r_{k}-r_{k+1}
$$
for some decreasing sequence $r_k>0$ starting at $r_0$. Geometrically, the second inequality implies that the map $I+P_k$ sends the ball $\{|z|\leq r_{k+1}\}$ into the ball $\{|z|\leq r_{k}\}$.
Assuming these inequalities, we deduce that 
$$\hat h_1(r_1)\leq (I+\hat P_0)(r_1)\leq r_1+(r_0-r_1)=r_0,$$
and then, by recurrence, that 
$$
\hat h_{k+1}(r_{k+1})\leq \hat h_k(r_{k+1}+\hat P_k(r_{k+1}))\leq \hat h_k(r_k)\leq r_0.
$$
Setting $r_{\infty}:=\lim r_k$, we deduce that 
$\hat h_k (r_{\infty})\leq r_0$ for all $k$, and since $h_k$ converges strongly to $h$, we deduce that 
$$\hat h(r_{\infty})\leq r_0,
$$
which implies that $\rho(h)\geq r_{\infty}$, and more precisely that the  map $h$ sends the ball 
$\{|z|\leq r_{\infty}\}$ into the ball $\{|z|\leq r_{0}\}$.

\subsubsection{}\label{cle}
Let us now enter the detail of the estimates and define the sequence  $r_k$. 
Reasoning by induction, we assume that $\hat F_k(r_k)\leq r_k$, and we will find an appropriate $r_{k+1}\in ]0, r_k[$ such that 
$$
\widehat{P_k}(r_k)\leq r_k-r_{k+1}, \quad \widehat{F_{k+1}}(r_{k+1})\leq r_{k+1}.
$$
For $i\leq 2^{k+1}$ we have 
$$|(P_k)_i|=|(F_k)_i|/|\lambda ^i-\lambda|= |(F_k)_i|/\omega_{i-1}\leq |(F_k)_i|/\Omega_{2^{k+1}}.
$$
Recalling that $P_k$ is a polynomial of degree $2^{k+1}$, we deduce that
$$\hat P_k(r)\leq \alpha_k^{-1} \hat F_k(r)
$$
where we set 
$$\alpha_k :=  \Omega_{2^{k+1}}.$$
Using that $F_k\in O(1+2^k)$, yields, for each $\gamma\in ]0,1[$, 
$$
\widehat{F_k}(r)\leq \gamma^{2^k}  \quad \forall r\in ]0, \gamma r_k]
$$
and therefore,
$$
\widehat{P_k}(r)\leq \gamma^{2k} \alpha_k^{-1} r,\quad \forall r\in ]0, \gamma r_k].
$$
We will have to apply this with a certain value $\gamma_k$ of the parameter $\gamma$. This value has to be chosen such that 
$\gamma^{2k} \alpha_k^{-1}<1$ by a certain margin. We  set 
$$
a_k:= \min (1/10, 1/k^2), \quad \gamma_k :=(\alpha_ka_k)^{2^{-k}}
$$
in such a way that 
$$
\widehat{F_k}(r)\leq a_k\alpha_k r, \quad \widehat{P_k}(r)\leq a_k r,  \quad \forall r\in ]0, \gamma_k r[.
$$
Note that  $\gamma_k$ depends strongly on the multiplier via the sequence $\alpha_k$. 
Using  paragraph \ref{rec}, we now  estimate the inverse $(I+R_k)=(I+P_k)^{-1}$:
$$
\hat R_k(r)\leq a_k(1-a_k)^{-1} r, \quad \forall r\in ]0, (1-a_k)\gamma_k r_k[.
$$
We can finally estimate 
\begin{align*}
	F_{k+1}=&(I+R_k)\circ(\lambda I+F_k)\circ (I+P_k)-\lambda I\\
	=&\lambda P_k+F_k\circ(I+P_k)+R_k\circ(\lambda I+F_k)\circ(I+P_k)
\end{align*}
by
$$
\widehat {F_{k+1}}(r)\leq \big(a_k+\alpha_k a_k(1+a_k)+a_k(1+a_k)(1+\alpha_ka_k)(1-a_k)^{-1}
\big)r \leq r
$$
(where the second inequality holds because $a_k\leq 1/10, \alpha_k\leq 2$)
for $r\leq r_{k+1}$ with 
$$ r_{k+1}:=(1-a_k)(1+\alpha_ka_k)^{-1}(1+a_k)^{-1}\gamma_k r_k.$$
With this choice of $r_{k+1}$, we just proved that $\widehat {F_{k+1}}(r_{k+1})\leq r_{k+1}$. Moreover
$$
\widehat{P_{k}}\leq a_k r_k \leq r_k -r_{k+1}
$$
where the last inequality holds because 
$r_{k+1}\leq (1-a_k)r_k$.
 This completes the proof of the key inequalities.

\subsubsection{}\label{arithm}
 We have proved that $\rho(h)\geq r_{\infty}$, it remains to describe how  $r_{\infty}$ depends on the 
multiplier $\lambda$, with (recalling that $a_k=\min (1/10, 1/k^2))$
\begin{align*}
r_{\infty}&=r_0\Pi_{k\geq 0} (1-a_k)(1+\Omega_{2^{k+1}}a_k)^{-1}(1+a_{k})^{-1}a_{k}^{2^{-k}}\Omega_{2^{k+1}}^{2^{-k}}\\
&\geq r_0\Big(\Pi_{k\geq 0}\Omega_{2^{k+1}}^{2^{-k}}\Big)
\Big(\Pi_{k\geq 0} (1-a_k)(1+2a_k)^{-1}(1+a_{k})^{-1}a_{k}^{2^{-k}}\Big)\\
&= Cr_0 \exp(-2b(\lambda)).
\end{align*}
Here $r_0$ is a scaling factor depending only  on the nonlinearity $F$, 
$$C:=  \Big(\Pi_{k\geq 0} (1-a_k)(1+2a_k)^{-1}(1+a_{k})^{-1}a_{k}^{2^{-k}}\Big)$$
is a universal positive constant,
and
$$b(\lambda) := \sum _{k\geq 1}2^{-k}\ln  \Omega_{2^k}^{-1}
$$
depends only on the multiplier $\lambda$ and is finite if and only if $\lambda$ is Bruno.
The constant $C$ is positive because the sums of general terms $-\ln (1-a_k)\sim a_k$, $\ln(1+a_k)\sim a_k$,
$2^{-k}\ln a_k $ are convergent. We see here that we had some freedom in the choice of the sequence $a_k$
(one could try to choose $a_k$ better in order to optimize the constant $C$).

\end{document}